\documentclass{amsart}

\input{epsf}
\usepackage{graphicx}
\usepackage{psfrag}
\usepackage{epsfig}

\newtheorem{theorem}{Theorem}[section]
\newtheorem{lemma}[theorem]{Lemma}

\theoremstyle{definition}
\newtheorem{definition}[theorem]{Definition}

\newtheorem{proposition}[theorem]{Proposition}
\newtheorem{conjecture}[theorem]{Conjecture}
\newtheorem{claim}[theorem]{Claim}

\theoremstyle{remark}
\newtheorem{remark}[theorem]{Remark}

\numberwithin{equation}{section}

\newcommand{\Z}{\mathbb{Z}}
\newcommand{\N}{\mathbb{N}} 
\newcommand{\K}{\mathcal{K}}

\newcommand{\B}{\mathcal{B}}
\newcommand{\s}{\sigma}
\newcommand{\de}{\delta }
\begin{document}
\title{Khovanov-Rozansky homology and the braid index of a knot}
\author{Keiko Kawamuro}
\address{Department of Mathematics, Rice University, Houston, Texas 77005}
\email{keiko.kawamuro@rice.edu}



\begin{abstract} We construct knots for whom the new Khovanov-Rozansky-Morton-Franks-Williams inequality gives a sharp bound for its braid index; however, the classical Morton-Franks-Williams inequality fails to do so. We also construct infinitely many knots for which the KR-MFW inequality fails to detect the braid indices.   \end{abstract}

\maketitle

\section{Introduction}
The Alexander's theorem states that any knot or link is isotopic to the closure of a braid. We can measure the complexity of knot $\K $ by the minimal possible number of braid strands, which is called the {\em braid index} $b_\K $. Morton \cite{Morton}, Franks and Williams \cite{FW} found an inequality which gives a lower bound for the braid index.

Let us first fix some notation.

Let $\K \subset S^3$ be an oriented knot or link and $\B_\K$ be the infinite set of closed braid diagrams of $\K.$ For an oriented closed braid diagram $D,$ let $b_D$ denote the number of the braid strands and $w_D$ the {\em writhe}, which is the number of positive crossings minus the number of negative crossings of $D$. Then the {\em braid index}  of $\K$ is $b_\K:=\min \{ b_D  | D \in \B_\K \}.$

We adopt the following definition and normalization of the HOMFLYPT polynomial $P_\K(a, q)$ defined by the skein relation; 
\begin{equation}\label{homfly-def}
a P_{\K_-}(a,q)  - a^{-1} P_{\K_+}(a,q) = (q- q^{-1}) P_{\K_0}(a,q) \quad  \mbox{ and } \quad P_{\mbox{unknot}}(a,q)=1. \end{equation} 
Let $d_\pm (\K)$ be the maximal (resp. minimal) $a$-degree of $P_\K(a, q).$ 

Now we state the Morton-Franks-Williams (MFW) inequality:

\begin{theorem}\label{MFW} {\em \cite{FW} \cite{Morton}}  {\em [\textbf{Morton-Franks-Williams inequality}] } For any closed braid diagram $D \in \B_\K$ of knot or link $\K$ we have $$w_D - b_D +1 \leq d_-( \K ) \leq d_+( \K ) \leq w_D + b_D -1.$$ Moreover $$\frac{d_+( \K ) - d_-( \K )}{2}+1\leq b_\K.$$  \end{theorem}

Khovanov-Rozansky homology \cite{KR} is a {\em categorification} of HOMFLYPT polynomial. 
In this paper, we use the {\em reduced} HOMFLY homology $\overline{H}^{i,j,k}( \K )$ of $\K$ introduced by Rasmussen \cite{R}. The graded Euler characteristic of the reduced HOMFLY homology is equal to the normalized HOMFLYPT polynomial (\ref{homfly-def});
$$\sum_{i,j,k} (-1)^{(k-j)/2}a^j q^i \dim \overline{H}^{i,j,k}( \K ) = P_\K (a,q).$$

Dunfield, Gukov, Rasmussen \cite{DGR} and Wu \cite{Wu} found Khovanov-Rozansky homology version of MFW inequality. We call it {\em KR-MFW inequality}.

\begin{theorem}\label{KR-MFW}{\em [ {\bf KR-MFW-inequality} ] \cite{DGR}, \cite{Wu}.}
Let $\K$ be a knot or a link and $$\de_+(\K) := \max \{ j \ | \ \overline{H}^{i,j,k}(\K) \neq 0, \mbox{ for some } i, k \},$$ $$\de_-(\K) := \min \{ j \ | \ \overline{H}^{i,j,k}(\K) \neq 0, \mbox{ for some } i, k \}.$$ Then, for any closed braid diagram $D \in {\B_\K}$ of $\K$ we have $$w_D - b_D +1 \leq \de_-(\K) \leq \de_+(\K) \leq w_D + b_D -1 \qquad \mbox{i.e., }  \ \ \frac{\de_+(\K)-\de_-(\K)}{2}+1\leq b_\K.$$
\end{theorem}

\begin{definition}\label{def-sharpness}
The MFW $($resp. KR-MFW$)$ inequality is called {\em sharp} on $\K$ if there exists $D \in {\B_\K}$ such that equalities $w_D - b_D +1 = d_-(\K)$ $($resp. $\de_-(\K))$ and $d_+(\K) = w_D + b_D -1$ $($resp. $\de_+(\K))$ hold. \end{definition}

Since $\de_-(\K) \leq d_-(\K) \leq d_+(\K) \leq \de_+(\K),$ we have the following.
\begin{proposition}
The sharpness of the MFW $($resp. KR-MFW$)$ inequality implies $$\frac{d_+(\K) -d_-(\K)}{2}+1 = b_\K \qquad (\mbox{resp. } \frac{\de_+(\K)-\de_-(\K)}{2}+1 = b_\K).$$
\end{proposition}

Elrifai \cite{E} has enumerated all the $3$-braids on which the MFW inequality is non-sharp.

\begin{theorem}\label{Elrifai} {\em \cite{E}  [{\bf Elrifai's example}] }
On all knots and links of braid index $=3$ the MFW inequality is sharp except 
\begin{eqnarray*}
\K_k &:=& \mbox{ the braid closure of }  (\s_1\ \s_2\ \s_2\ \s_1)^{2k}\  \s_1\ \s_2^{-2k-1}  \\
{\mathcal L}_k &:=& \mbox{ the braid closure of }  (\s_1\ \s_2\ \s_2\ \s_1)^{2k+1} \ \s_1\ \s_2^{-2k+1} \end{eqnarray*}
for $k\in \N$ and their mirror images $\overline{\K_k}, \ \overline{\mathcal L}_k.$
\end{theorem}

As the Euler characteristic of the KR homology gives us the HOMFLYPT polynomial, KR-homology contains more information than HOMFLYPT polynomial. It is interesting to find concrete examples that show the ``gap'' between KR-homology and HOMFLYPT polynomial. Elrifai's example seem to be natural candidates to see such gap. In fact, we have:

\begin{theorem}\label{main}
Let $$\K^\star := \K_1 = \mbox{ the braid closure of }  (\s_1\ \s_2\ \s_2\ \s_1)^2 \  \s_1\ \s_2^{-3}. $$
On $\K^\star$ and its mirror image $\overline{\K^\star}$ the MFW-inequality is not sharp but the KR-MFW inequality is sharp.
\end{theorem}

These are the first examples which show that Khovanov-Rozansky homology is ``stronger'' than HOMFLYPT polynomial in terms of detecting the braid index.

However, we will also see that Khovanov-Rozansky homology is not almighty. We study an obstruction of sharpness of the KR-MFW inequality and give infinitely many (and also first known) examples in Theorem~\ref{non-sharp-thm} whose braid index KR-homology fails to detect. 

Let $BM_{x,y,z,w}$ where $x, y, z, w \in \Z$ be the closure of the $4$-strand braid $$\sigma_1^x \ \sigma_2^y \ \sigma_3^{-1} \ \sigma_2^z \ \sigma_1^w \ \sigma_2 \ \sigma_3 \ \sigma_2 \ \sigma_2 \ \sigma_3.$$ It has been known \cite{K} that the set of BM-braids contains the five knots $9_{42},$ $9_{49},$ $10_{132},$ $10_{150},$ $10_{156},$ on which the MFW inequality is not sharp \cite{Jones}. Furthermore, the diagram contains infinitely many four tuples $(x,y,z,w)$ where the MFW inequality is not sharp  \cite{K}. A parallel result holds for the KR-MFW inequality:

\begin{theorem}\label{non-sharp-thm}
There are infinitely many four tuples $(x,y,z,w)$ such that the KR-MFW inequality is not sharp on $BM_{x,y,z,w}.$
\end{theorem}

Elrifai's examples have another interesting feature regarding to the generalized Jones' conjecture (\cite{MT}, \cite{K}) and the maximal Bennequin number Conjecture as we state below. See \cite{Jones} p.357 for Jones' original conjecture.

\begin{figure}[htpb!]
\begin{center}
\psfrag{g}{$w_\K+b_\K$}  
\psfrag{f}{$w_\K-b_\K$} 
\psfrag{c}{$w_\K$}  
\psfrag{b}{$b_\K$}
\psfrag{j}{writhe} 
\psfrag{i}{braid index}  
\psfrag{x}{$b$}  
\psfrag{y}{$w$}
\includegraphics[height=35mm]{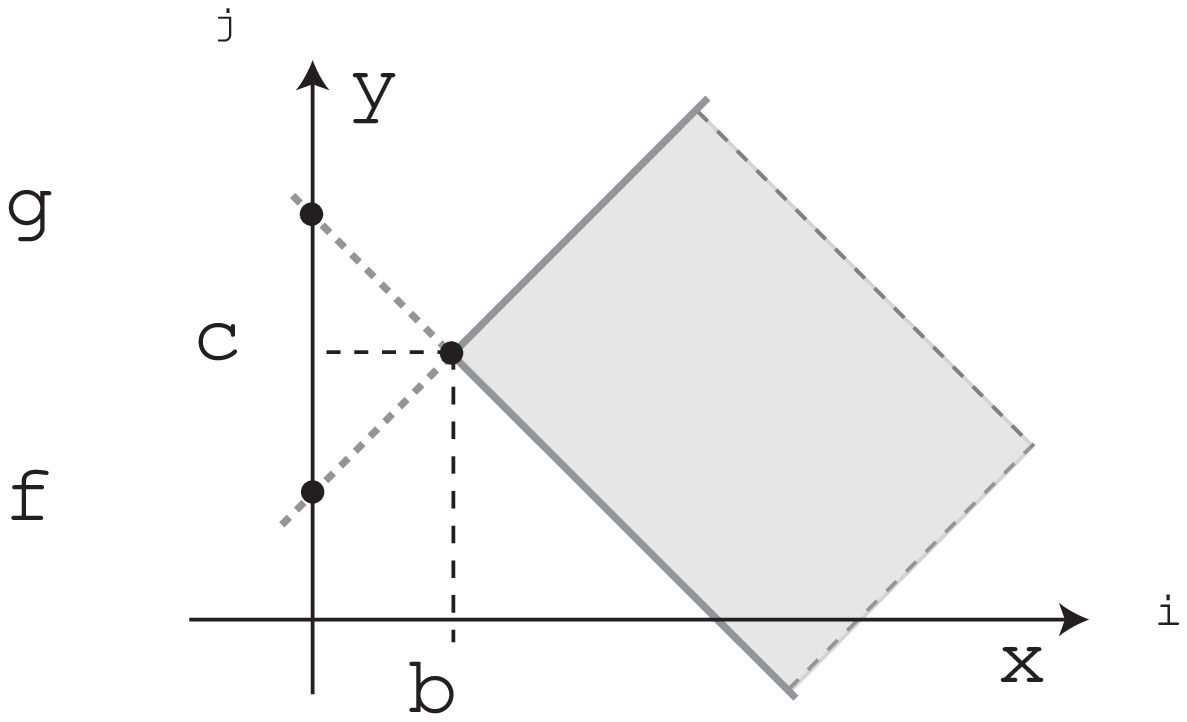}
\caption{The region of braid representatives of $\K$}\label{b-region}
\end{center}
\end{figure}

\begin{conjecture}\label{conj-k} [\textbf{generalized Jones' conjecture}] 
Let $\Phi:{\mathcal B}_{\mathcal K} \to \mathbb{N} \times \mathbb{Z}$ be a map with $\Phi(D):=\left( b_D, w_D \right)$ for $D\in {\mathcal B}_{\mathcal K}.$ Then there exists a unique $w_{\mathcal K} \in \mathbb{Z}$ such that
\begin{equation}\label{quadrant equation}
          \Phi({\mathcal B}_{\mathcal K})=\left\{ (b_{\mathcal K}+x+y, w_{\mathcal K}+x-y)\ |\ x,y\in \mathbb{N} \right\},
\end{equation}
which is the subset of the infinite shaded region shown in Figure \ref{b-region}.
\end{conjecture}

We can apply Conjecture \ref{conj-k} to contact geometry. Bennequin \cite{Ben} proved that any transversal knot in the standard contact structure $(S^3, \xi_{std})$ can be identified with a closed braid in $\mathbb{R}^3.$ The {\em Bennequin} ({\em self linking}) number $\beta_D := w_D - b_D$ of braid diagram $D\in \B_\K$ is an invariant of transversal knots and links in $(S^3, \xi_{std}).$ Then we denote the {\em maximal} Bennequin number of $\K$ by $\beta_\K := \{ \beta_D \ | \  D\in \B_\K \}.$ Conjecture \ref{conj-k} implies the following:

\begin{conjecture}\label{conj-k2} 
[\textbf{The maximal Bennequin number conjecture}] The maximal {\em Bennequin} number $\beta_\K$ of $\K$ is realized at the minimal braid index and $\beta_\K=w_\K - b_\K.$ 
\end{conjecture}

By Definition~\ref{def-sharpness}, it follows that:
\begin{proposition}\label{prop-sharp} 
The sharpness of the $($KR-$)$MFW inequality implies Conjectures~\ref{conj-k} and \ref{conj-k2}. 
\end{proposition}
Namely, the two inequalities $w_D - b_D +1 \leq d_-$ and $d_+ \leq w_D + b_D -1$ in Theorem~\ref{MFW} correspond to the two boundary lines of the infinite region in Figure~\ref{b-region}.

Conjecture \ref{conj-k} holds for many classes of knots and links: Franks and Williams \cite{FW} have proved the sharpness of the MFW inequality of knot and link with a braid representative of full positive twists with a positive word ($\Delta^{2n} P$) including unlinks, torus links and Lorenz links. Murasugi \cite{Murasugi} has affirmed the sharpness for alternating fibred knots and links and $2$-bridge knots and links. Jones  \cite{Jones} verified the sharpness for all the knots of less than or equal to $10$ crossing in the standard knot table except $9_{42}, 9_{49}, 10_{132}, 10_{150}, 10_{156}.$ Thus, Conjecture~\ref{conj-k} holds for these knots and links by Proposition~\ref{prop-sharp}.

In \cite{K2}, we have proved Conjecture \ref{conj-k} for $9_{42}, 9_{49}, 10_{132}, 10_{150}, 10_{156},$ by computing deficits of the MFW inequality of cabled knots,. 

In \cite{K}, more general results are given: If Conjecture~\ref{conj-k} holds for $\K, \mathcal{L}$ then Conjecture \ref{conj-k} also holds for the connect sum $\K \# {\mathcal L}$ and the $(p,q)$-cable of $\K.$

Thanks to Elrifai (Theorem~\ref{Elrifai}) we can improve the study of Conjectures~\ref{conj-k} and \ref{conj-k2} for the set of $3$-braids $B_3.$ Let $B_3' := B_3 \setminus \{ \ \K_k, {\mathcal L}_k, \overline{\K_k}, \overline{\mathcal L}_k \ | \ k \in \N \ \}.$ Rudolph's slice Bennequin inequality \cite{Rudolph} will play an important role to prove the following:
\begin{theorem}\label{prop-1}    
Conjecture~\ref{conj-k} holds for $B_3', \ \K_k$ and $\overline{\K_k}$ for $k\in \N.$ \end{theorem}
\begin{theorem}\label{prop-2}   
Conjecture~\ref{conj-k2} holds for $B_3', \ \K_k, {\mathcal L}_k$ and $\overline{\K_k}$ for $k\in \N.$ \end{theorem}

The rest of the paper is organized as follows. In Section~\ref{section 2}, we compute Khovanov-Rozansky homology of $\K^\star$ and prove Theorem~\ref{main}. In Section~\ref{section 3}, we discuss about non-sharpness of the KR-MFW inequality and prove Theorem~\ref{non-sharp-thm}. In Section~\ref{section 4}, we prove Theorems \ref{prop-1} and \ref{prop-2}.

\textbf{Acknowledgments } 
The author would like to thank Joan Birman, Michel Boileau, John Etnyre, Nathan Ryder and Hao Wu for helpful conversations and comments. She is especially grateful to Jacob Rasmussen for explaining KR-homology and pointing out Claim \ref{claim}.   


\section{Proof of Theorem~\ref{main}}\label{section 2}

We first recall some of the works of Rasmussen. Let $\sigma(\K)$ be the {\em signature} of $\K$. A knot $\K$ is called  {\em KR-thin} if $\overline{H}^{i,j,k}(\K)=0$ whenever $i+j+k\neq \sigma(\K).$  If $\mathcal{L}$ be a $2$-component link, we call it {\em KR-thin} if $\overline{\overline{H}}_N(\mathcal{L}),$ the {\em totally reduced} homology of ${\mathcal L}$, is thin for all sufficiently large $N>1.$ Denote $\de:=i+j+k$. In~\cite{R}, Rasmussen proved the following. 

\begin{theorem}\label{ras-work}
\begin{enumerate}
\item  Let $\K_+, \K_-, \K_0$ be links or knots differ by a single site. If $\K_-, \K_0$ are KR-thin and $$\det \K_- + 2 \det \K_0 = \det \K_+$$ then $\K_+$ is also KR-thin $($Corollary 7.7 of {\em \cite{R}).} 

\item  The connect sum of two KR-thin knots is also KR-thin $($Corollary 7.9 of {\em \cite{R}).}

\item   Among the knots with less than or equal to $9$ crossings, only $8_{19}, 9_{42}, 9_{43}, 9_{47}$ are {\em not} KR-thin $($Proposition 7.10 of {\em \cite{R}).}

\end{enumerate}
\end{theorem}

Now we prove Theorem~\ref{main}.

\noindent {\em Proof of Theorem~\ref{main}. }
We specify a resolution of $\K^\star.$ For simplicity let $n:=\sigma_n$ $(n=1,2)$ the generator of the Artin's braid group $B_3$ and $\overline{n}:=\sigma_n^{-1}.$ Let 
\begin{eqnarray*}
\K^\star= \K_+ &=&  \{1, 2, 2, 1, 1, 2, 2, 1, 1, \overline{2}, \overline{2}, \overline{2} \}      \\
\K_- &=&  \{1, 2, 2, 1, 1, 2, 2, 1, \overline{1}, \overline{2}, \overline{2}, \overline{2} \}  =  \{1, 2, 2, 1, 1, \overline{2} \} = \overline{5_2} \ \mbox{mirror image } \\
\K_0  &=&   \{1, 2, 2, 1, 1, 2, 2, 1, \overline{2}, \overline{2}, \overline{2} \} \\
\K_{0-} &=& \{1, 2, 2, 1, 1, 2, \overline{2}, 1, \overline{2}, \overline{2}, \overline{2} \}  = \{1, 2, 2, 1, 1, 1, \overline{2}, \overline{2}, \overline{2} \}\\
\K_{0 0} &=& \{1, 2, 2, 1, 1,2,1,\overline{2},\overline{2},\overline{2} \} =  \overline{5_2} \\  
\K_{0- -} &=& \{1, 2, \overline{2}, 1, 1, 1, \overline{2}, \overline{2}, \overline{2} \} =  \{1, 1, 1, 1, \overline{2}, \overline{2}, \overline{2} \} = T_{2, -3}  \# T_{2,4} \\ 
\K_{0- 0} &=& \{1, 2, 1, 1, 1, \overline{2}, \overline{2}, \overline{2} \}  = \mbox{unknot} 
\end{eqnarray*}
where $T_{p,q}$ is the $(p,q)$-torus knot or link.

\begin{figure}[htpb!]
\begin{center}
\psfrag{-}{$-$} \psfrag{0}{$0$} \psfrag{A}{$\K^\star= \K_+$} \psfrag{B}{$\K_-=\overline{5_2}$} \psfrag{c}{$\K_0  $} \psfrag{D}{$\K_{0-}$} \psfrag{E}{$\K_{0- -}= T_{2, -3}  \# T_{2,4} $} \psfrag{F}{$\K_{0- 0}=$unknot} \psfrag{G}{$\K_{00}=\overline{5_2}$} 
\includegraphics[height=28mm]{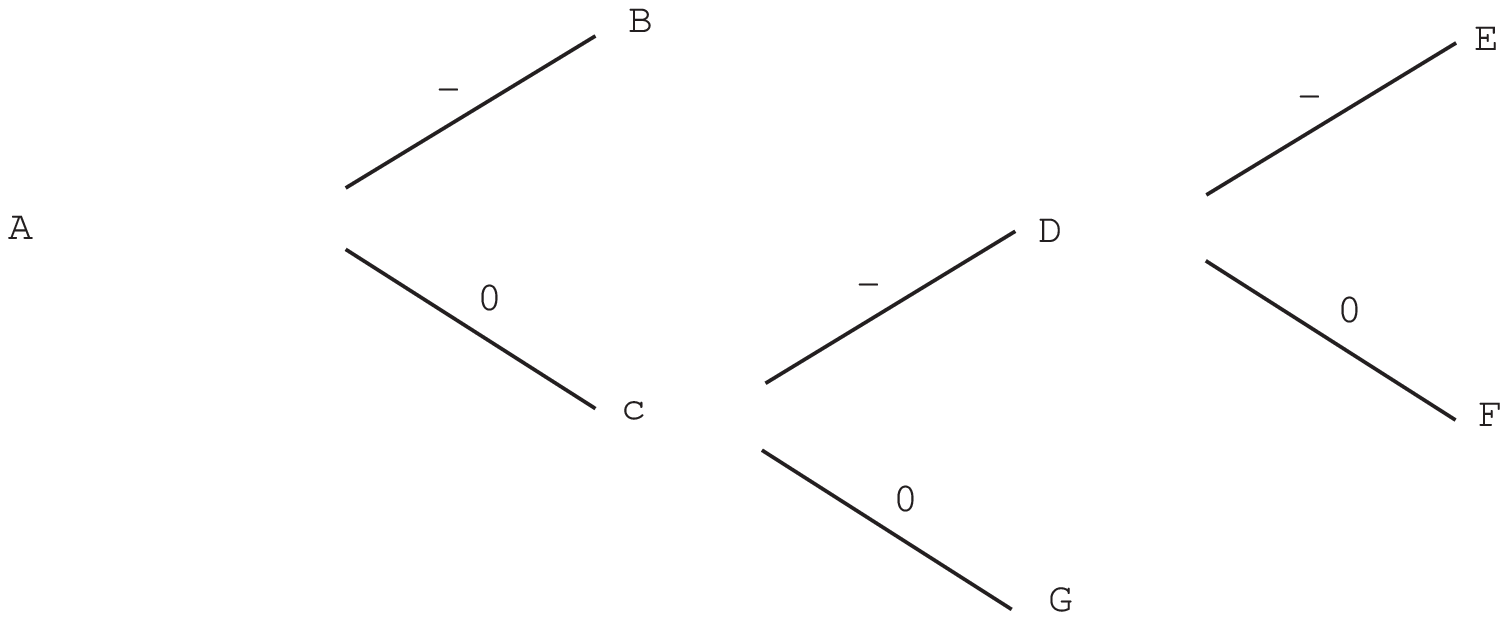}
\caption{Resolution of $\K^\star$ }\label{resolution}
\end{center}
\end{figure}
Thanks to Theorem~\ref{ras-work}.2 and \ref{ras-work}.3, knots $\overline{5_2}$ and $T_{2, -3}  \# T_{2,4}$ are KR-thin. Since $$\det(\K_{0--})+ 2 \det(\K_{0-0})=12 + 2 = 14=\det(\K_{0-}),$$ Theorem~\ref{ras-work}.1 tells that $\K_{0-}$ is also KR-thin and have the following table. 

\begin{center}
\begin{tabular}{lllccc}
\hline
      &  &  & $\sigma $  & $\det$  & KR-thin \\ \hline
$\K_+ $  & & knot & $2$  & $7$ & non-thin  \\ \hline
$\K_- = \K_{0 0}$ & $\overline{5_2}$ &  knot & $2$ & $7$ & thin  \\ \hline
$\K_0 $ & & link &  $1$ & $0$ & non-thin \\ \hline
$\K_{0-} $ & & link & $1$ & $14$ & thin \\ \hline
$\K_{0- -}$ & $T_{2, -3}  \# T_{2,4}$ & link & $1$ & $12$  & thin \\ \hline
$\K_{0- 0}$ & unknot & knot &  $0$  & $1$ & thin  \\  \hline
\end{tabular}
\end{center}

Let $\overline{H}_N(\K)$ be the reduced $sl(N)$ homology group defined by Khovanov and Rozansky \cite{KR}. It satisfies $$\sum_{I,J} (-1)^J q^I \dim \overline{H}_N^{I,J} (\K)=P_\K(q^N, q).$$ Let $M$ be a free module of rank $=4$ whose graded Poincare polynomial is $(q/t + t/q)^2.$ 
Using Proposition 7.6 of \cite{R} we have the following skein exact sequences.
\begin{equation}\label{ex sq 1}
\longrightarrow \overline{\overline{H}}_N (\K_{0-}) \longrightarrow \overline{H}_N (\K_{00}=\overline{5_2}) \otimes M \longrightarrow \overline{\overline{H}}_N (\K_{0}) \longrightarrow \overline{\overline{H}}_N (\K_{0-}) \longrightarrow
\end{equation}
\begin{equation}\label{ex sq 2}
\longrightarrow \overline{H}_N (\K_{-}=\overline{5_2}) \longrightarrow \overline{\overline{H}}_N (\K_{0}) \longrightarrow \overline{H}_N (\K_+) \longrightarrow \overline{H}_N (\K_{-}=\overline{5_2})  \longrightarrow
\end{equation}

Their HOMFLYPT polynomials satisfy:
\begin{eqnarray}
P_{\K_+}                   &=&  a^8 ( - q^4 -1 - q^{-4}) ) + a^6 (q^6  + q^2 + q^{-2} + q^{-6}) \label{homfly+}\\
P_{\overline{5_2}}                    &=&  -a^6 + a^4 (q^2 -1+q^{-2} )+ a^2 (q^2 -1+q^{-2} ) \notag \\
(q-q^{-1})^2 P_{\overline{5_2}}   &=& a^6 (-q^2 + 2 - q^{-2}) + a^4 ( q^4 - 3 q^2 + 4 - 3 q^{-2} + q^{-4}) \notag  \\ 
                                 && + a^2 ( q^4 - 3 q^2 + 4 - 3 q^{-2} + q^{-4}) \notag  \\
(q-q^{-1}) P_{\K_0}    &=&  a^7 ( q^4 + q^{-4} ) - a^5 (q^6 + 1 + q^{-6} ) + a^3 (q^2 -1 + q^{-2} ) \notag  \\
(q-q^{-1}) P_{\K_{0-}} &=& a^5(q^4 -q^2 +2-q^{-2} + q^{-4}) \notag \\
                                  & &+ a^3 (-q^6 + q^4 - 3 q^2 + 3 - 3q^{-2} + q^{-4} - q^{-6} )\notag  \\ 
                                 & & +  a (q^4 - 2q^2 + 3 -2q^{-2} + q^{-4} ) \notag 
\end{eqnarray}

\begin{figure}[htpb!]
\begin{center}
\psfrag{i}{$i$} \psfrag{j}{$j$} \psfrag{a}{$a$} \psfrag{b}{$b$} \psfrag{c}{$c$} \psfrag{d}{$d$} \psfrag{e}{$e$} \psfrag{f}{$f$} \psfrag{g}{$g$} \psfrag{h}{$h$} \psfrag{0}{$0$} \psfrag{2}{$2$} \psfrag{4}{$4$} \psfrag{6}{$6$} \psfrag{8}{$8$} \psfrag{-2}{$-2$} \psfrag{-4}{$-4$} \psfrag{-6}{$-6$} \psfrag{A}{$a'$} \psfrag{B}{$b'$} \psfrag{C}{$c'$} \psfrag{D}{$d'$} \psfrag{E}{$e'$}
\includegraphics [height=60mm]{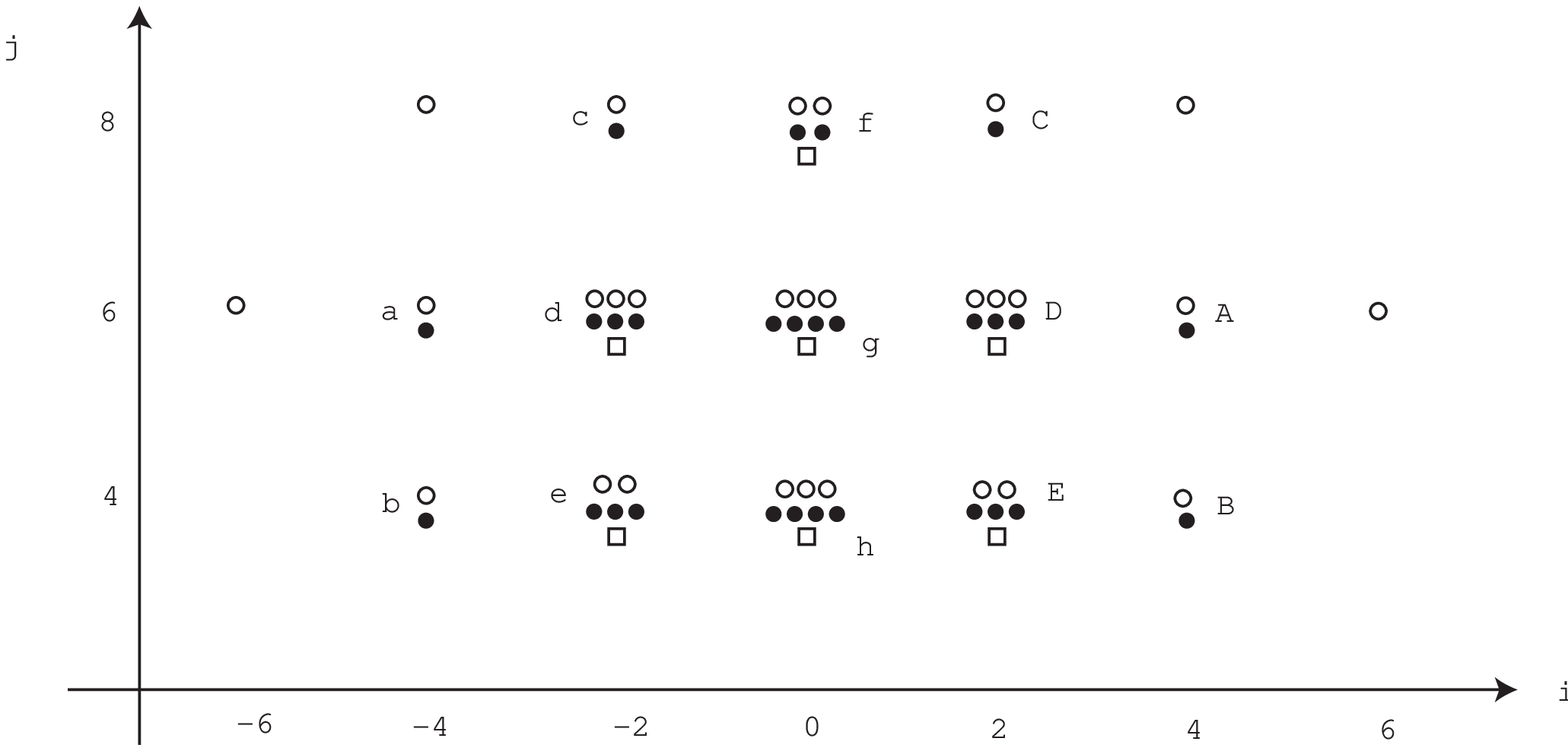}
\end{center}
\caption{Possible generators of $\overline{H}(\K^\star).$  The $\de$-gradings are $\de(\circ)=\de(\Box)=2, \ \de(\bullet)=4.$} \label{generators}
\end{figure}

Due to Rasmussen \cite{R}, there is a spectral sequence $E_k(N)$ whose $E_1$ term is $\overline{H}(\K)$ and converges to $\overline{H}_N(\K).$ When $N$ is large we have $\overline{H}(\K) \simeq \overline{H}_N(\K).$ Since $\overline{5_2}$ and $\K_{0 -}$ are KR-thin, we can explicitly compute $\overline{\overline{H}}_N (\K_{0-}), \ \overline{H}_N (\overline{5_2}) \otimes M$ and  $\overline{H}_N (\overline{5_2})$ from their HOMFLYPT polynomials. By the exact sequences (\ref{ex sq 1}) and (\ref{ex sq 2}), we guess possible generators of $\overline{H}(\K_+)$ with suitable $j$-grading shifts as illustrated in Figure \ref{generators}.  Hollow dots ($\circ$) represent $\overline{\overline{H}}(\K_{0-})$. Solid dots  ($\bullet$) represent $\overline{H}(\overline{5_2})\otimes M$. Squares ($\Box$) represent $\overline{H}(\overline{5_2})$.

We introduce the following claim, whose proof will be given shortly:

\begin{claim}\label{claim} The two generators at $b$ in Figure \ref{generators} survive. \end{claim}

Assuming Claim \ref{claim} and by the polynomial (\ref{homfly+}), we get $\de_-(\K^\star)=4$ and $\de_+(\K^\star)=8.$ Thus the KR-MFW inequality is sharp on $\K^\star.$
\hfill $\Box$

It remains to establish Claim~\ref{claim}.

\noindent {\em Proof of Claim~\ref{claim}. } 
According to \cite{S}, the Poincare polynomial of the reduced Khovanov homology $\overline{H}_{N=2}(\K^\star)$ is 
\begin{equation}\label{poly}
q^4+q^4t+q^6t^2+q^8t^2+q^8t^3+q^{10}t^3+2q^{10}t^4 + q^{12}t^5+q^{14}t^5+2q^{14}t^6
+ q^{16}t^7+q^{18}t^8+q^{20}t^9.
\end{equation}
Rasmussen \cite{R} proved that $$\overline{H}_N^{I,J}(\K) \simeq \bigoplus_{{i+Nj=I} \atop {(k-j)/2=J}} \overline{H}^{i,j,k}(\K)$$ where $I$ is the $q$-degree and $J$ the homological degree.
The only possible generators in Figure~\ref{generators} corresponding to the term $q^4 + q^4t$ of the polynomial (\ref{poly}) is the two at position $b$ since $I=4=-4+8=i+2j.$
\hfill $\Box$


\section{Proof of Theorem~\ref{non-sharp-thm}}\label{section 3}
\begin{lemma}\label{non-sharp-lemma}
Suppose that $D\in \B_\K$ is a closed braid diagram of ${\mathcal K}$ with $b_D = b_{\mathcal K}$. Focus on one site of $D$ and construct $D_+, D_-, D_0$ $($one of the three must be $D)$. Let $\alpha, \beta, \gamma \in \{ +, -, 0 \}$ be mutually distinct. Suppose $D_\alpha = D.$ 

If positive destabilization is applicable $p$-times to each of $ D_\beta, \ D_\gamma$, then 
\begin{equation}\label{1} (w_D + b_D -1) - \de_+(\K) \geq  2p. \end{equation}

If negative destabilization is applicable $n$-times to each of $ D_\beta, \ D_\gamma$, then 
\begin{equation}\label{2} \de_-(\K) - (w_D - b_D +1) \geq  2n. \end{equation}

Therefore, if $p+n>0$ the KR-MFW inequality is not sharp on ${\mathcal K}.$
\end{lemma}

\begin{proof}
Let $\K_\alpha, \K_\beta, \K_\gamma$ be the topological knot types of $D_\alpha, D_\beta, D_\gamma$ respectively. Thanks to Rasmussen's skein exact sequence of KR-homologies (Proposition 7.6 in \cite{R}), we have
$$\min \{\de_{-}(\K_{-})+2, \quad \de_{-}(\K_{0})+1\} \leq \de_{-}(\K_{+}) \leq \de_{+}(\K_{+}) \leq  \max \{\de_{+}(\K_{-})+2, \quad \de_{+}(\K_{0})+1\} $$
$$\min \{\de_{-}(\K_{+})-2, \quad \de_{-}(\K_{0})-1\} \leq \de_{-}(\K_{-})  \leq \de_{+}(\K_{-}) \leq  \max \{\de_{+}(\K_+)-2, \quad \de_{+}(\K_{0})-1\} $$
$$\min \{\de_{-}(\K_{+})-1, \quad \de_{-}(\K_{-})+1\} \leq \de_{-}(\K_{0}) \leq \de_{+}(\K_{0}) \leq  \max \{\de_{+}(\K_+) -1, \quad \de_{+}(\K_{-})+1\}. $$
Suppose $D_\alpha=D_+, \ D_\beta=D_-$ and $D_\gamma=D_0.$ Let $\tilde{D_-}$ (resp. $\tilde{D_0}$) be a closed braid diagram obtained by  applying positive destabilizations $p$ times to $D_-$ (resp. $D_0$). Then by the KR-MFW inequality we have
\begin{eqnarray*}
\de_+(\K_-) + 2 &\leq &( w_{\tilde{D}_-} + b_{\tilde{D}_-} -1 ) + 2   \\
&=& \{(w_{D_-} + b_{D_-} -1) -2p \} + 2  \\
&=&   (w_{D_+} -2) + b_{D_+} -1 -2p + 2   \\
&=&   (w_{D_+} + b_{D_+} -1) -2p, 
\end{eqnarray*}
and
\begin{eqnarray*}
\de_+(\K_0)+1 &\leq& (w_{\tilde{D_0}} + b_{\tilde{D_0}} -1 ) + 1 \\
&=&    (w_{D_0} + b_{D_0} -1 - 2p ) + 1\\
&=&    (w_{D_+} -1 ) + b_{D_+} -1 - 2p + 1 \\
&=&    (w_{D_+} + b_{D_+} -1) -2p.
\end{eqnarray*}
Thus,
$$\de_+ (\K_+)  \leq  \max \{\de_+(\K_-)+2, \quad \de_+(\K_0)+1\} \leq (w_{D_+} + b_{D_+} -1) -2p. $$
This is the inequality (\ref{1}) of the lemma. When $D_\alpha=D_-$ or $D_0,$ the same argument works.

Similarly, the inequality (\ref{2}) holds. 
\end{proof}

Now we are ready to prove Theorem \ref{non-sharp-thm}.

\noindent{\em Proof of Theorem~\ref{non-sharp-thm}. } Let $n:=\sigma_n$ and $\overline{n}:=\sigma_n^{-1} \in B_3$ where $n=1,2.$ By braid isotopy and destabilizations, define closed braid diagram  $ M_+,  M_-,  M_0$ by; 
\begin{eqnarray*}
 M_+ (\ =BM_{x,y,z,w} ) &:= &  2^x \ 3^y \ \overline{1}  \ \overline{2} \ 2^z \ 1^w \ 2 \ 3 \ 2 \ \overline{1},  \\ 
M_- & := &   2^x \ 3^y \ \overline{1}  \ \overline{2} \ 2^z \ 1^w \ 2 \ 3 \ 2 \ 1 \ \stackrel{+}{\Longrightarrow} \ 1^x \ 2^{y+1} \ 1^2 \ 2^{z+1} \ 1^w \ 2, \\ 
M_0 & := &   2^x \ 3^y \ \overline{1}  \ \overline{2} \ 2^z \ 1^w \ 2 \ 3 \ 2 \ \stackrel{+}{\Longrightarrow}  \ 2^y \ 1^{z+1} \ 2^{x+1} \ 1^{w+1} \ 2. 
\end{eqnarray*}

Here $\stackrel{+}{\Longrightarrow}$ means a positive destabilization and braid isotopy. Note that $b_{M_-}=b_{M_0}=3.$ 

In the proof of Theorem 2.8 in \cite{K}, it is proved that there exist infinitely many $(x,y,z,w)$'s such that the braid index of the topological knot type of $BM_{x,y,z,w}$ is $4.$ Thus Theorem~\ref{non-sharp-thm} follows from Lemma \ref{non-sharp-lemma}. \hfill $\Box$

\begin{remark}
As in Figure 14 of \cite{R}, Rasmussen explicitly computes the reduced KR-homology $\overline{H}(9_{42})$ of $9_{42}$ and we can see that the KR-MFW inequality is not sharp on $9_{42}.$ 
\end{remark}


\section{Maximal Bennequin numbers of $\K_k$ and ${\mathcal L}_k$}\label{section 4}

In this section, we prove Theorems~\ref{prop-1} and \ref{prop-2} together. The next lemma is essentially due to Elrifai~\cite{E}.

\begin{lemma}\label{lemma-Elrifai}
HOMFLYPT polynomial of $\K_k$ $($resp. ${\mathcal L}_k )$ coincides with the one for the $(2, 6k+1)$-torus knot $T_{2, 6k+1}$ $($resp. $T_{2, 6k+5}):$
$$P_{\K_k}(a, q)=P_{T_{2,6k+1}}(a, q), \qquad P_{{\mathcal L}_k}(a, q)=P_{T_{2, 6k+5}}(a, q). $$
\end{lemma}

\noindent{\em Proof of Theorems~\ref{prop-1} and \ref{prop-2}. } By Theorem~\ref{Elrifai} and Proposition~\ref{prop-sharp}, it follows that Conjectures~\ref{conj-k} and \ref{conj-k2} hold for the class $B_3'.$

Knots $\K_k, {\mathcal L}_k$ have the following quasipositive closed braid representations;
\begin{eqnarray*}
D_{\K_k} &:=& \overline{2}\ ( 1\ 2\ 2\ 1\ \overline{2} )^{2k}\ 1, \\
D_{{\mathcal L}_k} &:=& ( 1\ 2\ 2\ 1\ \overline{2} )^{2k-1}\ ( 1\ 2\ 2\ 1 )^2\ 1. 
\end{eqnarray*}
Since $(2\ 1\ \overline{2})$ and $(\overline{2}\ 1\ 2)$ are quasi positive factors, diagram $D_{\K_k} $ has $6k$ quasi positive factors and $D_{{\mathcal L}_k}$ has $6k+6.$

Let $g_4(\K)$ be the slice genus of $\K \subset S^3.$ If $\K$ has a quasipositive representative, say, $$D = ( w_1 \sigma_{j_1} w_1^{-1} ) \ (w_2 \sigma_{j_2} w_2^{-1} ) \cdots (w_p \sigma_{j_p} w_p^{-1} ),$$ then  thanks to Rudolph \cite{Rudolph} we have $2\ g_4(\K) = p-b_D + 1.$ Therefore,
\begin{eqnarray*}
2\ g_4(\K_k) &=& 6k-3+1=6k-2, \\
2\ g_4({\mathcal L}_k) &=& (6k+6)-3+1=6k+4.
\end{eqnarray*}

Next recall Rudolph's {\em slice-Bennequin Inequality} \cite{Rudolph};
$$ 2\ g_4(\K) \geq w_D - b_D + 1 = \beta_D +1\quad \mbox{ for } \quad D \in \B_\K.$$
On each $D_{\K_k},\  D_{{\mathcal L}_k},$ the inequality is sharp;
\begin{eqnarray*}
2\ g_4(\K_k) = 6k-2 = w_{D_{\K_k}} - b_{D_{\K_k}} + 1, \\
2\ g_4({\mathcal L}_k) = 6k+4 = w_{D_{{\mathcal L}_k}} - b_{D_{{\mathcal L}_k}} + 1.
\end{eqnarray*}

Since we know that $D_{\K_k},  D_{{\mathcal L}_k}$ have $b_{\K_k}=b_{D_{\K_k}}=3$ and $b_{{\mathcal L}_k}=b_{D_{{\mathcal L}_k}}=3$ (the minimal braid index), Conjecture \ref{conj-k2} holds for $\K_k$ and ${\mathcal L}_k.$

By Lemma~\ref{lemma-Elrifai}, we can compute the minimal and the maximal $a$-degrees of $P_{\K_k}(a, q)$ and $P_{{\mathcal L}_k}(a, q):$
$$d_-(\K_k)=6k, \quad d_+(\K_k)=6k+2,$$
$$d_-({\mathcal L}_k)=6k+4, \quad d_+({\mathcal L}_k)=6k+6.$$
Since $w_{D_{\K_k}} + b_{D_{\K_k}} -1=6k+2,$ one of the MFW inequalities $d_+ \leq w_D + b_D -1$ is sharp on $\K_k.$ This sharpness combined with the above argument about the slice Bennequin inequality, we conclude that Conjectures~\ref{conj-k} and \ref{conj-k2} hold for $\K_k$ and their mirror image $\overline{\K_k}.$

We remark that since $w_{D_{{\mathcal L}_k}} + b_{D_{{\mathcal L}_k}} -1=6k+8 > 6k+6 =d_+({\mathcal L}_k)$  (the MFW inequality is not sharp on ${\mathcal L}_k),$ the same argument does not apply to ${\mathcal L}_k.$
\hfill $\Box$

\bibliographystyle{amsplain}

\end{document}